\documentclass[11pt,a4paper,leqno]{amsart}
\usepackage{epsfig,graphicx}
\newtheorem{thm}{Theorem}[section]

\newtheorem{prop}[thm]{Proposition}
\newcommand{\e}{\epsilon}
\newcommand{\de}{\delta}

\begin{document}

\title{$\theta$-curve polynomials and finite-type invariants}
\author[Y. Huh]{Youngsik Huh}
\author[G.T. Jin]{Gyo Taek Jin}
\address{Multimedia Lab, SAIT, P.O.\ Box 111, Suwon 440-600, Korea}
\email{yshuh@samsung.com}
\address{Department of Mathematics, KAIST, Taejon 305-701, Korea}
\email{trefoil@kaist.ac.kr}

\maketitle
\thispagestyle{empty}
\bigskip\bigskip

\begin{abstract}
The normalized Yamada polynomial, $\widetilde{R}_A $, is a
polynomial invariant in variable A for $\theta$-curves. In this
work, we show that the coefficients of $\widetilde{R} _{e^{x}}$
which is obtained by replacing A with $e^x=\sum x^n/n!$ are
finite-type invariants for $\theta$-curves although the
coefficients of original $\widetilde{R}_A$ are not finite-type. A
similar result can be obtained in the case of Yokota polynomial
for $\theta$-curves.
\end{abstract}

\section{Introduction}
Birman and Lin discovered infinitely many finite type invariants
for knots derived from polynomial\footnote{In this article, a
polynomial stands for a Laurent polynomial.} invariants~
\cite{BL}. They showed that the coefficients of one variable
HOMFLY polynomial and Kauffman polynomial are finite-type with
substituting the variable $t$ by $e^x$. Bar-Natan showed that the
coefficients of Conway polynomial are also finite-type \cite{BN}.
On the other hand, Zhu observed that the coefficients of Jones
polynomial are not finite-type \cite{ZHU}. Using Zhu's idea, Jin
and Lee showed that the coefficients of the 2-variable HOMFLY
polynomial, the 2-variable Kauffman polynomial and the
Q-polynomial are not finite type invariants~\cite{JL}.

In this paper, following the studies listed above, we examine
whether the coefficients of two polynomial invariants of
$\theta$-curves are finite-type invariants. A $\theta$-curve is a
graph embbedded in $\mathbf R^3$ consisting of two vertices and
three edges between them.

Originally, finite-type invariants were defined for knots. Later,
Stanford extended them to links and some other spatial graphs
including $\theta$-curves and Kanenobu investigated them with
emphasis on $\theta$-curves \cite{S1,S2,K}. In fact, abundant
finite-type invariants for $\theta$-curves can be obtained in the
following way: There is a well-defined 3-component link
associated with each $\theta$-curve, which is the boundary of a
surface obtained by thickening the $\theta$-curve in a canonical
way~\cite{KSWZ}. Every finite-type invariant for this associated
link was shown to be a finite type invariant for the original
$\theta$-curve~\cite{S2}. Therefore all the finite type
invariants derived from the polynomial invariants of the
3-component links associated with $\theta$-curves are finite type
invariants.

On the other hand, there are some invariants for $\theta$-curves
which seem to have different origins from the above. One of them
is Yamada polynomial~\cite{YA}. There is no known relation between
the Yamada polynomial of a $\theta$-curve and invariants of its
associated 3-component link. Yamada polynomial can be calculated
by skein relations from diagrams while this method doesn't seem to
work for the invariants of $\theta$-curves coming from their
associated 3-component links.

Let $\widetilde{R}_A$ be the Yamada polynomial for
$\theta$-curves normalized to behave multiplicatively under
connected sums. We will show that the coefficients of
$\widetilde{R}_{e^x}$ which is obtained by replacing the variable
$A$ with $e^x=\sum_{n=0}^\infty x^n/n!$ are finite-type
invariants for $\theta$-curves although the coefficients of the
original $\widetilde{R}_A$ are not finite-type. In 1999, Yokota
introduced a polynomial invariant for
$\theta_m$-curves~\cite{YO}. A $\theta_m$-curve is a spatial
graph consisting of 2 vertices and $m$-edges between them. This
polynomial is obtained from the $SU(m)$ invariant of a linear
combination of links constructed from a given $\theta$-curve
diagram. This method is different from that of~\cite {KSWZ}. For
the coefficients of Yokota polynomial, we obtain the same results
as in the case of Yamada polynomial.

\section{Definitions}
\begin{figure}[b]
\centering
\includegraphics[scale=0.7]{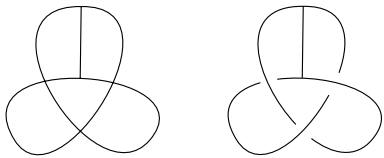}
\caption{Projection and diagram of $3_1$}\label{fig:3_1} %
\end{figure}

\begin{figure}[t]
\centering
\includegraphics[scale=0.85]{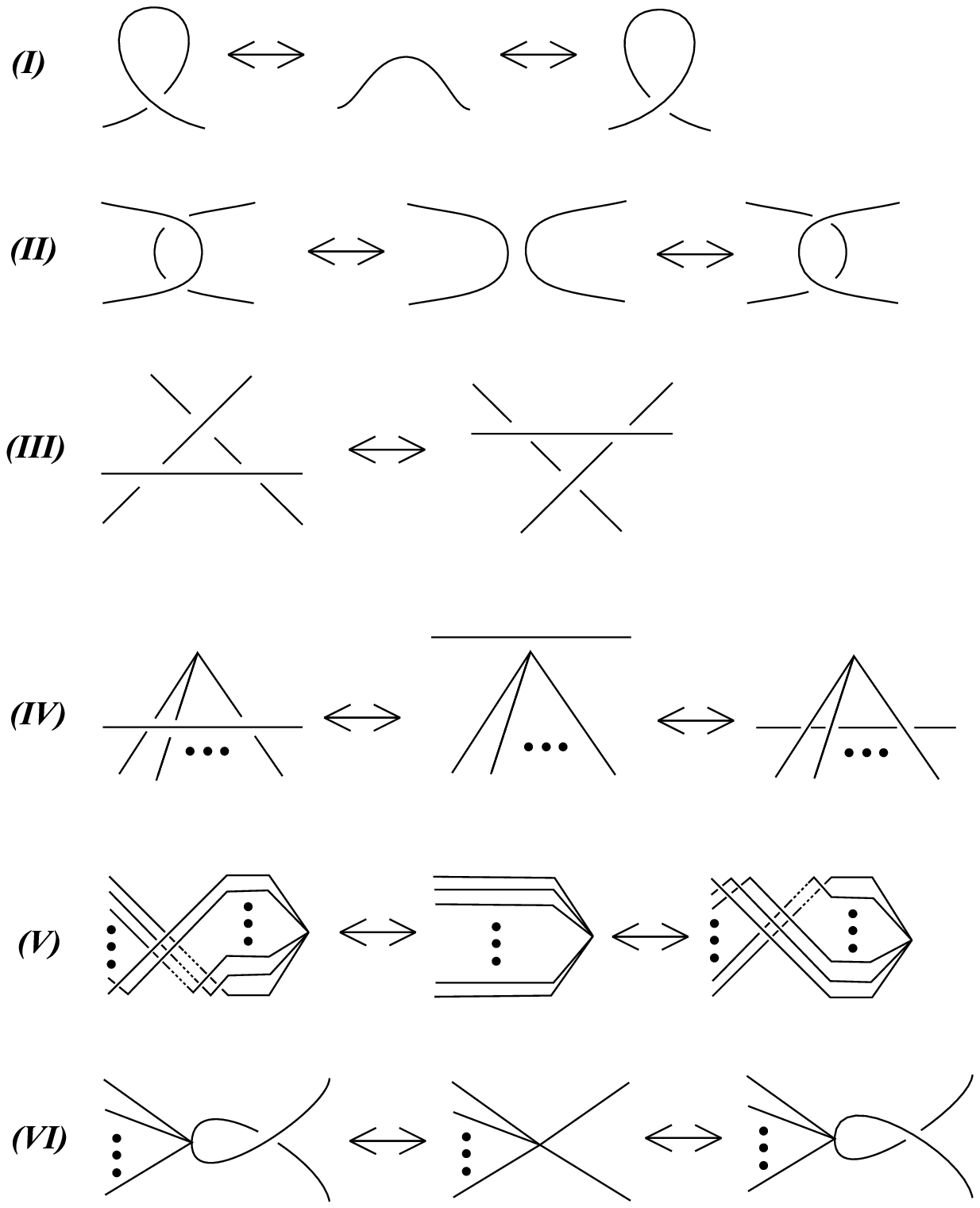} 
\caption{Moves}\label{fig:moves}
\end{figure}

A \emph{spatial graph\/} is a finite graph which is PL embedded in
${\mathbf R}^{3}$. Two spatial graphs $G$ and $G'$ are said to be
\emph{equivalent\/} or \emph{isotopic\/} if there exists an
orientation preserving auto\-homeomorphism of ${\mathbf R}^3$
carrying one onto to the other. A \emph{projection\/} of a spatial
graph is its image under a natural projection map of $\mathbf
R^3$ onto a euclidean plane. We always assume that spatial graphs
are in general position with the projection maps so that the only
singular points of the projections are transverse double points
away from vertices. A double point in a projection is called a
\emph{crossing}. A \emph{diagram\/} of a spatial graph is its
projection with informations denoting  which strand is over or
under at each crossing. See Figure~\ref{fig:3_1}. Two spatial
graphs are equivalent if and only if a diagram of one of them can
be transformed to a diagram of the other by a finite sequence of
moves I--VI given in Figure~\ref{fig:moves}~\cite{K1}.

\begin{figure}[h]
\centering
\includegraphics[scale=0.6]{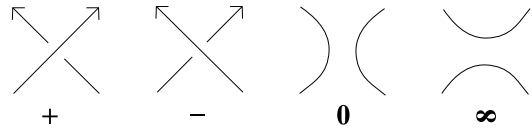} 
\caption{Positive, negative crossings and $0$, $\infty$ changes}
\label{fig:crossings}
\end{figure}

A \emph{$\theta$-curve\/} is a spatial graph consisting of two
vertices and 3 edges joining them. A $\theta$-curve is said to be
\emph{trivial\/} if it is equivalent to a $\theta$-curve in
${\mathbf R}^2$. A crossing in a diagram of a $\theta$-curve is
called a \emph{positive crossing\/} (resp.\ \emph{negative
crossing\/}) or said to have \emph{crossing number\/} `$+1$'
(resp.\ `$-1$') if it is like the one in
Figure~\ref{fig:crossings} with `$+$' sign (resp.\ `$-$' sign)
when all the edges are oriented coherently from one vertex to the
other.

For a positive integer $n$, an \emph{$n$-sign\/} is an $n$-tuple
of $1$, $-1$, $0$ and $\infty$, i.e., an element
$(\e_1,\ldots,\e_n)\in\{1,-1,0,\infty\}^n$. Let
$C=\{c_1,c_2,\ldots,c_n\}$ be a set of $n$ distinct crossings in a
diagram $D$ of a $\theta$-curve and let $\e =
(\e_1,\ldots,\e_n)$ be an $n$-sign. Denote by
$D_\e$ the diagram obtained from $D$ by replacing each
crossing $c_i$ with a local diagram of Figure~\ref{fig:crossings}
corresponding to $\e_i$.

For any rational invariant $v$ of $\theta$-curves, we define
 $$v(D\mid C)=\sum_{\e \in \{1,-1\}^n}(-1)^{|\e|} v(D_
 \e)$$
where $|\e|=\sum_{\e_i>0}\e_i$. We say that $v$
is {\em finite-type of order less than $n$} if $v(D\mid C)=0$ for
every set $C$ of $n$ distinct crossings in every $\theta$-curve
diagram~$D$. If this vanishing condition is true for any $n+1$
crossings but not for some $n$ crossings, we say that $v$ is
\emph{finite-type of order $n$}.

Given two $\theta$-curves $\Theta$ and $\Theta'$, we may assume
that $\Theta$ is in the upper half space ${\mathbf R}^3_+=\{
(x,y,z) \in {\mathbf R}^3\mid z \ge 0 \}$ except for a small
neighborhood of one vertex, $\Theta'$ is in the lower half space
${\mathbf R}^3_-=\{ (x,y,z) \in {\mathbf R}^3 \mid z \le 0 \}$
except for a small neighborhood of one vertex, $\Theta \cap
{\mathbf R}^2_0 = \Theta' \cap {\mathbf R}^2_0 = \{\mbox{ $3$
points } \}$ where ${\mathbf R}^2_0 = \{ (x,y,z) \in {\mathbf
R}^3 \mid z = 0 \}$, and $(\Theta\cap\mathbf
R_-^3)\cup(\Theta'\cap\mathbf R_+^3)$ is a trivial
$\theta$-curve. Then ($\Theta \cap {\mathbf R}^3_+$) $\cup$
($\Theta' \cap {\mathbf R}^3_-$) is a $\theta$-curve which is
called a \emph{connected sum\/} of $\Theta$ and $\Theta'$,
denoted by $\Theta \# \Theta'$.

\section{Yamada Polynomial}
In \cite{YA}, Yamada introduced a  polynomial invariant in
variable $A$ for diagrams of spatial graphs, which is invariant
under the moves II--IV in Figure~\ref{fig:moves} and is a flat
isotopy\footnote{For $\theta$-curves, flat isotopy implies
isotopy, since the the move VI at a trivalent vertex can be
derived from the moves I--V.} invariant up to a multiplication by
a power of $-A$. Two diagrams of spatial graphs are said to be
\emph{flat isotopic\/} if one can be transformed into the other
by a finite sequence of moves I--V in Figure~\ref{fig:moves}. For
a diagram $D$ of a $\theta$-curve, let $R_A(D)$ be the polynomial
invariant of Yamada for $D$. Then $R_A$ can be determined by the
moves II--IV and  the formulae Ya1--Ya5, where $\ominus$,
$\bigcirc\kern-1.5pt{-}\kern-1.5pt\bigcirc$ and `$D\ \bigcirc$'
denote a diagram of a $\theta$-curve in a plane, that of a
handcuff curve\footnote{A handcuff curve is a connected spatial
graph with two vertices and three edges obtained by connecting
two disjoint loop edges with the third edge.} in a plane and the
diagram $D$ with a disjoint unknotted circle, respectively, and
$\sigma=A+1+A^{-1}$.
\def\mybox#1{\parbox{0.6cm}{\epsfxsize=0.6cm \epsfbox{#1}}}
\begin{align}
\tag{Ya1} &R_A(\ominus)=\sigma-\sigma^2=-(A^2+A+2+A^{-1}+A^{-2}) \\ %
\tag{Ya2} &R_A\left(\mybox{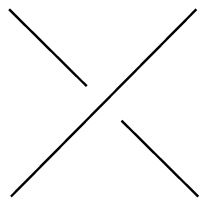}\right) -
R_A\left(\mybox{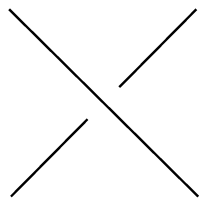}\right) = (A-A^{-1})\left(
R_A\left(\mybox{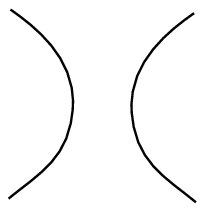}\right)
- R_A\left(\mybox{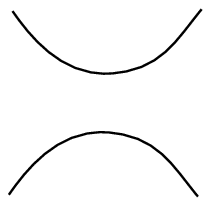}\right) \right) \label{YP:skein}\\ %
\tag{Ya3}
&R_A\left(\mybox{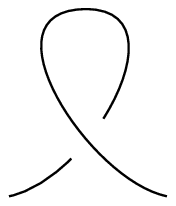}\right)=A^2R_A\left(\mybox{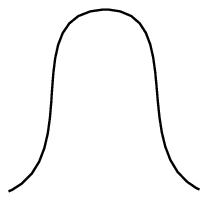}\right),\qquad
R_A\left(\mybox{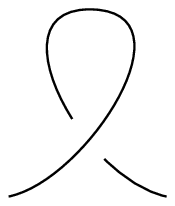}\right)=A^{-2}R_A\left(\mybox{R1_no.eps}\right) \\ %
\tag{Ya4}
&R_A\left(\mybox{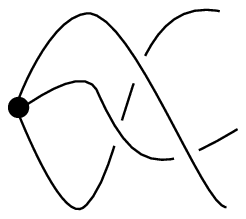}\right)=-A^{-3}R_A\left(\mybox{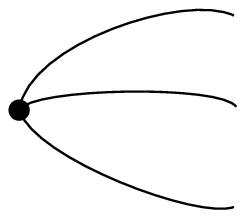}\right),\qquad
R_A\left(\mybox{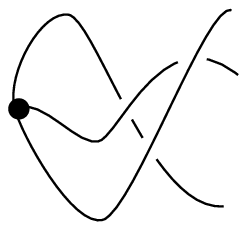}\right)=-A^{3}R_A\left(\mybox{T_no.eps}\right) \\ %
\tag{Ya5} &R_A(\bigcirc\kern-1.5pt{-}\kern-1.5pt\bigcirc)=0 \\ %
\tag{Ya6} &R_A(D\mbox\ \bigcirc)=\sigma R_A(D)
\end{align}

 A crossing in a diagram is called a \emph{self-crossing\/}
if its two strands are from a single edge, and a {\em
non-self-crossing\/} otherwise. For a diagram $D$ of a
$\theta$-curve, let $s(D)$ and $n(D)$  denote the sums of the
crossing numbers of the self-crossings in $D$ and the
non-self-crossings in $D$, respectively.

\begin{prop}\label{prop:R-tilde}
For any $\theta$-curve diagram $D$, we define
$$\widetilde{R}_A(D)=(-A)^{n(D)-2s(D)}R_A(D)/(\sigma-\sigma^2).$$ %
Then $\widetilde R_A$ is an isotopy invariant of $\theta$-curves
with the following properties:
\begin{enumerate}
\item \label{item:poly-A} $\widetilde R_A(D)$ is a polynomial in variable $A$.
\item \label{item:R-sum} $\widetilde{R} _A(\Theta \# \Theta')= \widetilde{R} _A(\Theta) \widetilde{R} _A(\Theta')$
for two $\theta$-curves $\Theta$, $\Theta'$.
\item \label{item:mirror}If $\bar D $ is the mirror image of $D$, then $\widetilde{R}_A(\bar D) = \widetilde{R}
_{A^{-1}}(D)$.
\end{enumerate}
\end{prop}
\begin{proof} According to \cite[Theorem~7]{YA}, the polynomial
$(-A)^{n-2s}R_A$ is an isotopy invariant of $\theta$-curves.
Therefore so is $\widetilde R_A$. The property
`\ref{item:poly-A}' is a consequence of the fact that $R_A$ is a
multiple of $\sigma-\sigma^2$ for any $\theta$-curve diagram $D$,
which is not hard to see.
The property `\ref{item:R-sum}' follows from the connected sum
formula
 $$R_A(D\sharp D^\prime)=R_A(D)R_A(D^\prime)/(\sigma-\sigma^2)$$
for any $\theta$-curve diagrams $D$ and
$D^\prime$~\cite[Theorem~5]{YA}.
Finally, the property `\ref{item:mirror}' is a consequence of the
fact $R_A(\bar D) = R_{A^{-1}}(D)$~\cite[Proposition~6]{YA}.
\end{proof}

\begin{thm}
\label{th1} For any integer $n$, the coefficient of $x^n$ in the
power series $\widetilde{R} _{e^x}$, obtained from $\widetilde{R}
_A$ by the substitution $A=e^x= \sum _{i=0} ^{\infty}x^i/i!$, is a
finite-type invariant of order at most $n$. On the other hand,
the coefficients of the polynomial $\widetilde{R} _A$, in variable
$A$, are not finite-type invariants.
\end{thm}
\begin{proof}
Let $D$ be a $\theta$-curve diagram. Suppose $D_+$, $D_-$, $D_0$
and $D_\infty$ are the diagrams with one crossing of $D$, say $c$,
replaced by the local diagrams of Figure~\ref{fig:crossings},
respectively. Then we have
\begin{align*}
 \widetilde{R}_A(D_+)&=(-A)^{m  }R_A(D_+)/(\sigma-\sigma^2)\\
 \widetilde{R}_A(D_-)&=(-A)^{m+j}R_A(D_-)/(\sigma-\sigma^2),
\end{align*}
where $m=m(D_+)=n(D_+)-2s(D_+)$ and
 $$j=j(c)=\begin{cases}
     4&\ \mbox{if $c$ is a self-crossing,}\\
     -2&\ \mbox{if $c$ is a non-self-crossing.}
 \end{cases}$$
In the following computations, we omit the subscript $A$ from $\widetilde R_A$ and $R_A$.

Using (\ref{YP:skein}), we obtain
\begin{align}
\tag{E1} \label{E:R-tilde}
\begin{split}
\lefteqn{ (\sigma - \sigma ^2)(\widetilde{R}(D_+)-
\widetilde{R}(D_-))} \\
 & =   (-A)^{m}(R(D_+ )-A^{j} R (D_-))  \\
 & =  (-A)^{m}((1-A^{j})R(D_-)  + (A-A^{-1})(R(D_0)-R(D_{\infty}))).
\end{split}
\end{align}

Suppose $D_+$, $D_-$, $D_0$ and $D_\infty$ are the diagrams
obtained from $D$ by a repeated application of the skein relation
(Ya2), which are identical except at one place where they differ
as indicated in Figure~\ref{fig:crossings}.
The local orientations of $D_+$ and $D_-$ shown in
Figure~\ref{fig:crossings} are induced from the orientation of
$D$ in which all edges are oriented from one vertex to the other.
Using
(\ref{YP:skein}) again, we have
\begin{align}
\tag{E2} \label{E:R}
\begin{split}
\lefteqn{R(D_+) - A^t R(D_-)}  \\
 & =  R(D_-) +(A-A^{-1})(R(D_0)-R(D_{\infty}))-A^tR(D_-) \\
 & =  (1- A^t)R(D_-)+(A-A^{-1})(R(D_0) -R(D_{\infty})),
\end{split}
\end{align}
for any nonzero integer $t$.

\medskip
For a set $C=\{c_1,c_2,\ldots,c_n\}$ of $n$ distinct crossings of
$D$, we consider the sum
\begin{equation}\tag{E3}\label{E:R(D|C)}
(\sigma-\sigma^2)\widetilde{R}(D\mid C)
  = \sum_{\e
  \in\{1,-1\}^n}(-1)^{|\e|}(\sigma-\sigma^2)\widetilde{R}(D_{\e}).
\end{equation}

For $1\le k\le n$, let $E_k=\{1,-1\}^k$ and $F_k=\{-1,0,\infty\}^k$.
For a $k$-sign $\e=(\e_1,\ldots,\e_k)\in\{1,-1,0,\infty\}^k$,
let $\e+$, $\e-$, $\e0$ and $\e\infty$
be the $(k+1)$-signs obtained  from $\e$ by appending $1$, $-1$, $0$ and $\infty$,
respectively.
Let $f_k=-(1-A^{j(c_k)})$ and $f_0=-f_\infty=-(A-A^{-1})$.
For each $k$-sign $\de=(\de_1,\ldots,\de_k)\in F_k$, with $1\le k\le n$, we define
$$G_i(\de)=\begin{cases}f_{n+1-i}&\mbox{if\ }\de_i=-1\\
                        f_0      &\mbox{if\ }\de_i=0\\
                        f_\infty &\mbox{if\ }\de_i=\infty,
           \end{cases}$$
for $1\le i\le k$. Applying (\ref{E:R-tilde}) to the right hand side of
(\ref{E:R(D|C)}), we obtain
\begin{align*} \tag{E4} \label{E:En}
\begin{split}
\lefteqn{\sum_{\e\in E_n}(-1)^{|\e|}(\sigma-\sigma^2)\widetilde{R}(D_{\e})}\\
 & = \sum_{\e\in E_{n-1}}(-1)^{|\e+|}(\sigma-\sigma^2)(\widetilde{R}(D_{\e+})-\widetilde{R}(D_{\e-}))   \\
 & = \sum_{\e\in E_{n-1}}(-1)^{|\e|}(-A)^{m(D_{\e+})}(f_n R(D_{\e-}) + f_0 R(D_{\e0}) +f_\infty R(D_{\e\infty}))\\
 &=\sum_{\de\in F_1}G_1(\de)\sum_{\e\in E_{n-1}}(-1)^{|\e|}(-A)^{m(D_{\e+})}R(D_{\e\de}).
\end{split} \end{align*}

\begin{figure}[b]
\includegraphics[scale=0.7]{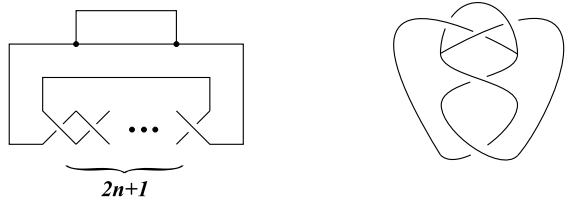} \caption{$T_{2n+1}$
and $5_1$}\label{fig:curves}
\end{figure}

If $\e\in E_{n-k-1}$,  then
\begin{align*}
& |\e{+}|=|\e{-}|+1=|\e|+1,\\
& m(D_{\e-+\cdots+})=m(D_{\e++\cdots+})+j(c_{n-k}),
\end{align*}
where `$+\cdots+$' means a duplication of $+$'s so that the length of the whole `sign' is $n$.
Using (\ref{E:R}) with a $k$-sign $\de\in F_k$, we ontain
\begin{align*}
\begin{split}
\lefteqn{\sum_{\e\in E_{n-k}}(-1)^{|\e|}(-A)^{m(D_{\e+\cdots+})}R(D_{\e\de})}\\
 & = \sum_{\e\in E_{n-k-1}}(-1)^{|\e+|}(-A)^{m(D_{\e+\cdots+})}(R(D_{\e{+}\de})-A^{j(c_{n-k})}R(D_{\e{-}\de}))   \\
 & = \sum_{\e\in E_{n-k-1}}(-1)^{|\e|}(-A)^{m(D_{\e+\cdots+})}(f_{n-k} R(D_{\e{-}\de})+f_0R(D_{\e{0}\de})+f_\infty R(D_{\e{\infty}\de}))\\
 &=\sum_{\de^\prime\in F_1}G_{k+1}(\de^\prime)\sum_{\e\in E_{n-k-1}}(-1)^{|\e|}(-A)^{m(D_{\e+\cdots+})}R(D_{\e\de^\prime\de}).
\end{split} \end{align*}
Repeated applications of this to (\ref{E:En}) lead us to
\begin{align*}
\begin{split}
\lefteqn{\sum_{\e\in E_n}(-1)^{|\e|}(\sigma-\sigma^2)\widetilde{R}(D_{\e})} \\
&=\sum_{\de\in F_2}G_1(\de)G_{2}(\de)\sum_{\e\in E_{n-2}}(-1)^{|\e|}(-A)^{m(D_{\e++})}R(D_{\e\de}) \\
&= \cdots \\
&=\sum_{\de\in F_{n-1}}\prod_{i=1}^{n-1}G_i(\de)\sum_{\e\in E_1}(-1)^{|\e|}(-A)^{m(D_{\e+\cdots+})}R(D_{\e\de}) \\
&=(-A)^{m(D_{+\cdots+})}\sum_{\de\in F_n}\prod_{i=1}^nG_i(\de)R(D_\de).
\end{split} \end{align*}
This shows that $(\sigma-\sigma^2)\widetilde{R}(D\mid C)$ is divisible by $(1-A)^n$,
since each $G_i(\de)$ is divisible by $1-A$.
Therefore, after the substitution $A=e^x=\sum_{i=0}^\infty
x^i/i!$, the power series $\widetilde{R}_{e^x}(D\mid C)$  has no
terms of degree less than $n$. This proves the first part of the
theorem.

Let  $v_r(\Theta)$ be the coefficient of $A^r$ in $\widetilde{R} (\Theta)$
for a $\theta$-curve $\Theta$. To show that $v_r$ is
not finite-type for any integer $r$, we use the $\theta$-curves
$T_{2n+1}$ with $n\ge0$ and $5_1$ shown in
Figure~\ref{fig:curves}.

Since $T_1$ is trivial, we have $\widetilde R(T_1)=1$.
Inductively, we can compute the maximal and the minimal degrees of
$\widetilde{R}(T_{2n+1})$. For $n \ge 1$, they are $-2n$ and
$-(8n+3)$, respectively.
For $n\ge1$, let $C_{2n}$ be a set of $2n$ distinct crossings of
$T_{4n+1}$. The following computation shows that $v_0$ is not
finite-type.
 $$v_0(T_{4n+1}\mid C_{2n}) = \sum_{p=0} ^{2n} (-1)^p \binom{2n}{p} v_0(T_{4n-2p+1}) =1 \neq 0.$$

Let $G = \bar T_5 \# 5_1$. Since the minimal degree\footnote
{$\widetilde R(5_1)=
A^{10}-2A^9-2A^8+6A^7-2A^6-6A^5+8A^4+A^3-7A^2+3A+3-3A^{-1}+A^{-3}$.}
of $\widetilde R(5_1)$ is $-3$, that of
$\widetilde{R}(G)=\widetilde{R}(\bar
T_5)\widetilde{R}(5_1)=\widetilde{R}_{A^{-1}}(T_5)\widetilde{R}(5_1)$
is~$1$.
For a positive integer $r$, let $G^r$ be a connect sum of $r$
copies of $ G$. Since the minimal degree of
 $$\widetilde{R}( G ^r \# \bar T_{2n+1})
= (\widetilde{R}( G ))^r \widetilde{R}(\bar T_{2n+1}) $$ is
$r+2n$, we see that $v_r( G ^r \#  \bar T_{2n+1}) \ne 0$ if and
only if $n=0$. This leads us to
 $$
 v_r( G ^r \# \bar T_{4n+1}\mid{C_{2n}})
 =\sum_{p=0}^{2n}(-1)^p\binom{2n}p v_r(G^r\#\bar T_{4n-2p+1})\ne0,$$
for any $n\ge1$. Considering the mirror images and the maximal
degrees, we also obtain
 $$
 v_{-r}(\bar G ^r \#T_{4n+1}\mid{C_{2n}})
 =\sum_{p=0}^{2n}(-1)^p\binom{2n}p v_{-r}(\bar G^r\#T_{4n-2p+1})\ne0,$$
for any $n\ge1$. This proves the second part of the theorem.
\end{proof}

\section{Yokota polynomial for $\theta _m$-curves}
In \cite{YO}, Yokota introduced a polynomial invariant for
$\theta_m$-curves. It is a normalization of the $SU(m)$ invariant
of a linear combination of link diagrams derived from a given
$\theta _m$-curve diagram $D$ in a way which is different from
that of~\cite {KSWZ}. The polynomial has also some properties
about local changes in diagram which enable us to compute it.
Among them, we introduce only what is necessary for
$\theta$-curves. Let $\langle D\rangle$ be the $SU(3)$ invariant
derived from a $\theta$-curve diagram $D$. We assume that the
edges of $D$ are oriented coherently from one vertex to the other.
Then $\left<D\right>$ is a polynomial in variable $t$ which is
invariant under the moves II--IV and is determined by the
following formulae:
\begin{align}
\tag{Yo1} &\left<\ominus\right>=1 \\ %
\tag{Yo2}
 & t\left<\mybox{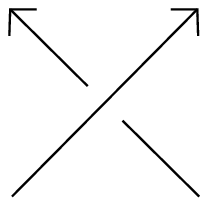}\right>-t^{-1}\left<\mybox{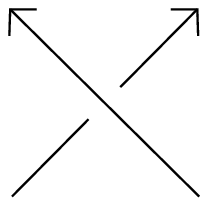}\right>
  =(t^3-t^{-3})\left<\mybox{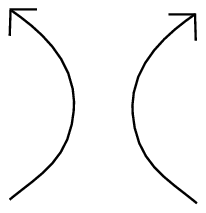}\right>\\
\tag{Yo3}
 &t^{-8}\left<\mybox{R1_po.eps}\right>=\left<\mybox{R1_no.eps}\right>=t^{8}\left<\mybox{R1_ne.eps}\right> \\ %
\tag{Yo4}
 &-t^{4}\left<\mybox{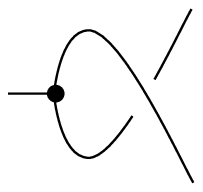}\right>=\left<\mybox{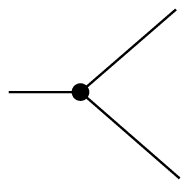}\right>=-t^{-4}\left<\mybox{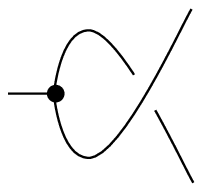}\right>\\ %
\tag{Yo5} &\left<D\ \bigcirc\right>=(t^6+1+t^{-6}) \left<D\right> %
\end{align}
\begin{prop}
\cite {YO} \label{yoko} For any $\theta$-curve diagram $D$,  we
define
 $$P(D) =
(-t^{4})^{n(D)-2s(D)}\langle D\rangle ,$$ Then, $P_z$, obtained
from $P$ by the substitution $z=t^3$, is an isotopy invariant for
$\theta$-curves with the following properties:
\begin{enumerate}
\item $P_z(D)$ is a polynomial in variable $z$.
\item $P_z(\Theta  \# \Theta' )=P_z(\Theta ) P_z(\Theta')$
for two $\theta $-curves $\Theta $ and $\Theta'$.
\item If $\bar D$ is the mirror image  of  $D$, then $P_z(\bar
D) = P_{z^{-1}}(D)$.
\end{enumerate}\end{prop}

\begin{thm}
\label{th2} For any integer $n$, the coefficient of $x^n$ in the
power series $P_{e^x}$, obtained from $P_z$ by the substitution
$z=e^x=\sum_{i=0}^\infty x^i/i!$, is a finite-type invariant of
order at most $n$. On the other hand, the non-trivial
coefficients, i.e., those of even degree terms, of the polynomial
$P_z$ are not finite-type invariants.
\end{thm}
\begin{proof} The first part can be proven in a  similar way to that of
Theorem~\ref{th1}, using the skein relation (Yo2). It is easily
seen that $P_z$ has only even degree terms. The maximal degree and
the minimal degree of $P_z(T_{2n+1})$ are $-4n$ and $-(8n+4)$,
respectively, for $n \ge1$, and
those\footnote{$P_z(3_1)=z^2+z^8-z^{10}$.} of $P_z(3_1)$ are $10$
and $2$, respectively. The second part can be proven in a similar
way to that of Theorem~\ref{th1}, using $3_1 ^r \#\bar T_{4n+1}$
instead of $ G ^r \# \bar T_{4n+1}$.
\end{proof}

\section*{Acknowledgement}
The first author would like to thank Yoshiyuki Yokota for helpful
answers to his questions.


\begin{thebibliography}{XXXX}
\bibitem[B] {BN} D. Bar-Natan, {\em On the Vassiliev knot
invariants}, Topology {\bf34} (1995) 423--472.
\bibitem[BL] {BL} J. S. Birman and X. S. Lin, {Knot polynomials
and Vassiliev's invariants}, Invent. Math. {\bf111} (1993) 225--270.
\bibitem[HJO1]{HJO1} Y. Huh, G.T. Jin and S. Oh, {\em Strongly almost
            trivial $\theta$-curves}, preprint.
\bibitem[HJO2]{HJO2} Y. Huh, G.T. Jin and S. Oh, {\em Elementary set
            for $\theta_n$-curve projections}, preprint.
\bibitem[JL] {JL} G.T. Jin and J. Lee, {\em Coefficients of
HOMFLY polynomial and Kauffman polynomial are not finite-type
invariants}, preprint.
\bibitem[Kan] {K} T. Kanenobu, {\em Vassiliev-type invariants of a
theta-curve}, J. Knot Theory Ramifications {\bf6}(4) (1997) 455--477.
\bibitem[Kau] {K1} L. H. Kauffman, {\em Invariants of graphs in
three-space}, Trans. Amer. Math. Soc. {\bf311} (1989) 697--710.
\bibitem[KSWZ] {KSWZ} L. H. Kauffman, J. Simon, K. Wolcott and
P. Zhao, {\em Invariants of theta-curves and other graphs in
3-space}, Topology Appl. {\bf49} (1993) 193--216.
\bibitem[S1] {S1} T. Stanford, {\em Finite-type invariants of
knots, links, and graphs}, Topology, {\bf35} (1996) 1027--1050.
\bibitem[S2] {S2} T. Stanford, {\em The functionality of Vassiliev-type
invariants of links, braids, and knotted graphs},
J. Knot Theory Ramifications {\bf3} (1994) 247--262.
\bibitem[W]{W} K. Wolcott, {\em The knotting of theta-curves and
            other graphs in $S^3$},
            Geometry and Topology, Marcel Decker (1987) 325--346.
\bibitem[Ya] {YA} S. Yamada, {\em An invariant of spatial graphs},
J. Graph Theory {\bf13} (1989) 537--551.
\bibitem[Yo]{YO} Y. Yokota, {\em Polynomial invariants of
            $\theta_m$-curves in 3-space}, a lecture at 7th Japan-Korea
            School of Knots and Links, Kobe, 1999.
\bibitem[Zha] {ZHA} P. Zhao, {\em Is knotted graph determined by its
associated links?}, Topology Appl. {\bf57} (1994) 23--30.
\bibitem[Zhu] {ZHU} J. Zhu, {\em On Jones knot invariants and
Vassiliev invariants}, New Zealand J. of Math. {\bf27} (1998),
293--299.
\end{thebibliography}
\end{document}